\documentclass[a4paper,10pt]{amsart}
\usepackage{t1enc,amssymb,epsfig}
\usepackage[english]{babel}

\input xy
\xyoption{all}

\newtheorem{proposition}{Proposition}[section]
\newtheorem{theorem}[proposition]{Theorem}
\newtheorem{lemma}[proposition]{Lemma}

\newtheorem{corollary}[proposition]{Corollary}
\newtheorem{remark}[proposition]{Remark}

\begin{document}
\title[The Nash map for surface singularities]{\textbf{A  class of
    non-rational surface
    singularities \\with bijective Nash map}}
\author{\sc Camille Pl{\'e}nat}
\address{Universit{\'e} de Provence, LATP UMR 6632\\
Centre de Math{\'e}matiques et Informatique, 39 rue Joliot-Curie,13453
    Marseille cedex 13, France.} 
\email{plenat@cmi.univ-mrs.fr}
\author{\sc Patrick Popescu-Pampu}
\address{Univ. Paris 7 Denis Diderot, Inst. de
  Maths.-UMR CNRS 7586, {\'e}quipe "G{\'e}om{\'e}trie et dynamique" \\case
  7012, 2, place Jussieu, 75251 Paris cedex 05, France.}
\email{ppopescu@math.jussieu.fr}


\thispagestyle{empty}
\begin{abstract}{Let $(\mathcal{S},0)$ be a germ of complex analytic
    normal surface. On its minimal resolution,
    we consider the reduced exceptional divisor $E$ and its
    irreducible components $E_{i}, \: i \in I$. The Nash
    map associates to each irreducible component $C_k$ of the space of
    arcs through $0$ on $\mathcal{S}$ the unique component of $E$ cut
    by the strict transform of the generic arc in $C_k$.
    Nash proved its injectivity and asked if it was bijective. As a
    particular case of our main theorem, we prove that this is the case if
    $E.E_{i} <0$ for any $i \in I$.}
\end{abstract}

\maketitle

\par\medskip\centerline{\rule{2cm}{0.2mm}}\medskip
\setcounter{section}{0}

\section{Introduction}

Let $(\mathcal{S},0)$ be a germ of complex analytic normal
surface. Let
$$\pi_m : (\tilde{\mathcal{S}}_m ,E) \rightarrow (\mathcal{S},0)$$
be its
minimal resolution, where $E$ is the reduced exceptional divisor of
$\pi_m$, and let $(E_i)_{i \in I}$ be the irreducible components of
$E$. A resolution is called \textit{good} if $E$ has normal
crossings and if all its components are smooth. It is important to
notice that, by Grauert's contractibility theorem for surfaces \cite{G
  62}, there exist
singularities whose minimal resolution is not good.

An \textit{arc through } $0$ on $\mathcal{S}$ is a germ of formal map
$(\mathbf{C},0) \rightarrow (\mathcal{S},0)$.

We denote by $(\mathcal{S},0)_{\infty}$ the \textit{space of
  arcs} through $0$ on $\mathcal{S}$. It can be canonically given the
structure of a scheme over $\mathbf{C}$, as the projective limit of
schemes of finite type obtained by truncating arcs at each finite
order. So it makes sense to speak about its irreducible components
$(C_k)_{k \in K}$. For each arc represented by an element in
  $C_k$, one can consider the intersection point with $E$ of its
strict transform on $\tilde{\mathcal{S}}_m$. For a generic element of $C_k$ (in
  the Zariski topology), this
  point is situated on a unique irreducible component of $E$. In
  this manner one defines a map:
    $$ \mathcal{N} : \{ C_k\: |  \: k \in K\} \rightarrow \{E_i \: |
    \: i \in I\}$$
which is called \textit{the Nash map} associated to the germ
    $(\mathcal{S},0)$. It was defined by Nash around  1966, in a preprint
    published later as \cite{N 95}. He proved that the map
    $\mathcal{N}$ is injective (which shows in particular that $K$ is
    a finite set) and asked the question:
\medskip

 \textit{Is the map $\mathcal{N}$
    bijective?}

\medskip

This question is now called \textit{the Nash problem on arcs}. No germ
$(\mathcal{S},0)$ is known for which the answer is negative. But the
bijectivity of $\mathcal{N}$ was only proved till now for special
classes of  singularities:

$\bullet$ for the germs of type $(\mathbf{A}_n)_{n \geq 1}$ by Nash himself in
\cite{N 95};

$\bullet$ for normal minimal singularities by Reguera \cite{R 95}; 
different proofs were given by Pl{\'e}nat \cite{P 04} and by
Fern{\'a}ndez-S{\'a}nchez \cite{F 05};

$\bullet$ for sandwiched singularities it was sketched by Reguera
\cite{R 04}, using her common work  \cite{LR 99} with
Lejeune-Jalabert on the wedge problem.

 $\bullet$ for the germs of type $(\mathbf{D}_n)_{n \geq 4}$ by
Pl{\'e}nat
 \cite{P  04};

$\bullet$ for the germs with a good $\mathbf{C}^*$-action such that
the curve $\mathrm{Proj} \: \mathcal{S}$ is not rational, it follows
immediately by combining results of Lejeune-Jalabert \cite{L 80} and
Reguera \cite{R 04}.

With the exception of the last class, all the other ones consist only
in \textit{rational} singularities and can be defined purely
topologically.
\medskip

Here we prove that the Nash map is bijective for a new class of
surface singularities (Theorem \ref{extension}), whose definition depends
only on the intersection matrix of the minimal resolution and not on
the genera or possible singularities of the components $E_i$. In
particular, their minimal resolution may not be good, which contrasts
with the classes of singularities described before. The
following (Corollary \ref{mthm2}) is a particular
case of the main theorem:

\medskip
 \textit{If $E.E_i <0$ for any $i \in I$, then the Nash map
 $\mathcal{N}$ is bijective.}
\medskip

We also show (Corollary \ref{infinity}) that the hypothesis of the
previous corollary are verified by an
infinity of topologically pairwise distinct \textit{non-rational}
singularities, which explains the title of the article.

The Nash map can also be defined in higher dimensions, over any field,
for not necessarily normal schemes which admit
resolutions of their singularities. It is always
injective and the same question can also be asked. Ishii and Koll{\'a}r
proved in \cite{IK 03} (a good source for everything we use about
spaces of arcs and the Nash map, as well as for references on related
works) that it is not always bijective. Indeed, they
gave a counterexample in dimension $4$, which can be immediately
transformed in a counterexample in any larger dimension. They left
open the cases of dimensions 2 and 3...

\subsection*{Acknowledgements}
   We are grateful to the organizers of the Third Franco-Japa\-nese
   Symposium and School on Singularities, Hokkaido University
   (Sapporo), September 13 - 17, 2004, during which the present
   collaboration was started. We are also grateful to
   M. Lejeune-Jalabert, M. Morales and the referee for their remarks. 

\medskip

\section{A criterion for distinguishing components of the space of
  arcs} \label{criterion}

\medskip

Consider a germ $(\mathcal{S},0)$ of normal surface and its minimal
resolution morphism $\pi_m : (\tilde{\mathcal{S}}_m, E)\rightarrow
(\mathcal{S},0)$. If $D$ is a divisor on $\tilde{\mathcal{S}}_m$, it can be
uniquely written as the
sum of a divisor supported by $E$ - called \textit{the exceptional
  part} of $D$ - and a divisor whose support meets $E$ in a finite
number of points. If $D$ consists only of its exceptional part, we say
that $D$ is \textit{purely exceptional}.

For each $i \in I$, let $v_{E_i}$ be the divisorial valuation defined
by $E_i$ on the fraction field of the analytic local ring
$\mathcal{O}_{\mathcal{S},0}$. Denote by
$\frak{m}_{\mathcal{S},0}$ the maximal ideal of this local ring. If
$f \in \frak{m}_{\mathcal{S},0}$, the exceptional part of
$\mathrm{div}(f \circ \pi_m)$ is precisely $\sum_{i\in I}
v_{E_i}(f)E_i$.

For each component $E_i$ of $E$, consider the arcs on $\tilde{\mathcal{S}}_m$
whose closed points are on $E_i- \cup_{j \neq i}E_j$ and which
intersect $E_i$ transversally. Consider the set of their images in
$(\mathcal{S},0)_{\infty}$ and denote its closure by
$V(E_i)$. The sets $V(E_i)$ are irreducible 
and $$(S,0)_\infty = \bigcup_{i \in I}V(E_i)$$ (see Lejeune-Jalabert
\cite[Appendix $3$]{LJ 90}).

The following proposition is a special case of a general one
proved by the first author in \cite{P 03} (see also \cite{P 04})
for non-necessarily normal germs of any dimension and for
arbitrary resolutions. It generalizes an equivalent result proved by
Reguera \cite[Theorem 1.10]{R 95} for the case of rational surface
singularities. 

\begin{proposition} \label{distinguish}
   If there exists a function $f \in \frak{m}_{\mathcal{S},0}$ such
   that \linebreak $v_{E_i}(f) < v_{E_j}(f)$, then $V(E_i) \nsubseteq
   V(E_j)$.
\end{proposition}

\textit{Proof:} Let $(\mathcal{S},0)\hookrightarrow (\mathbf{C}^n,0)$
be an analytic embedding of the germ  $(\mathcal{S},0)$. Denote by
$(x_1,...,x_n)$ the coordinates of $\mathbf{C}^n$. An arc $\phi \in
(\mathcal{S},0)_{\infty}$ is then represented by $n$ formal power
series $(x_k(t)= \sum_{l=1}^{\infty}a_{kl}t^l)_{1 \leq k \leq n}$,
where the coefficients $(a_{kl})_{k,l}$ are subjected to algebraic
constraints, coming from the fact that the arc must lie on
$\mathcal{S}$.

For each $j \in I$, a Zariski open set $U_{f}(E_j)$ in $V(E_j)$ 
consists of the images by
$\pi_m$ of the arcs on $\tilde{\mathcal{S}}_m$ which meet transversely
$E_j$ in a 
smooth point of $\mathrm{div}(f \circ \pi_m)$. If $\phi \in
U_{f}(E_j)$, we have:
  $$v_{E_j}(f) = v_t(f \circ \phi)$$
where $v_t(g)$ denotes the order in $t$ of $g \in \mathbf{C}[[t]]$.

This shows that the first $(v_{E_j}(f)-1)$ coefficients of $f\circ
\phi$, seen as elements of $\mathbf{C}[a_{kl}]_{k,l}$, must
vanish. Their vanishing defines a closed subscheme
$Z_{f,j}$ of $(\mathcal{S},0)_{\infty}$. Therefore, $U_{f}(E_j) \subset
Z_{f,j}$, which implies that:
   $$ V(E_j) \subset Z_{f,j}.$$

As $v_{E_i}(f) < v_{E_j}(f)$, we see that no element of
$U_{f}(E_i)$ is included in $Z_{f,j}$, which shows that:
   $$ V(E_i) \nsubseteq Z_{f,j}.$$

The proposition follows. \hfill $\Box$

\medskip

\section{Construction of functions with prescribed divisor} \label{prescribed}

In this section,  $\pi: (\tilde{\mathcal{S}}, E) \rightarrow
(\mathcal{S},0)$ denotes \textit{any} resolution of
$(\mathcal{S},0)$.

Inside the free abelian group generated by $(E_i)_{i
  \in I}$ we consider the set:
  $$\mathcal{L}(\pi):= \{D \: | \: D\neq 0, \: D\cdot E_{i}\leq 0, \:
  \forall \: i\in I\}.$$
It is a semigroup with respect to addition, which we call (following
  L{\^e} \cite[3.2.5]{L 98})  \textit{the Lipman
  semigroup} associated to
$\pi$ (see Lipman \cite[\S 18]{L 69}). It is known that it consists only
  of effective divisors (see Lipman \cite[\S 18 (ii)]{L 69}).

We call \textit{strict Lipman semigroup} of $\pi$ the subset:
$$\mathcal{L}^{\circ}(\pi):=\{D\in \mathcal{L}(\pi)\: | \: D\cdot E_i <0,
\: \forall \: i\in I\}$$
of the Lipman semigroup of $\pi$. It is always non-empty.

The importance of the Lipman semigroup comes from the fact that the
exceptional parts of the divisors of the form $\mathrm{div}(f \circ
\pi)$, where $f \in \frak{m}_{\mathcal{S},0}$, are elements of
it. The converse is true for rational surface singularities, but this is 
not the case for arbitrary surface singularities.  

We 
give now a numerical criterion on a divisor $D \in \mathcal{L}(\pi)$
which allows one to conclude that it is the exceptional part of a
divisor of the form $\mathrm{div}(f \circ \pi)$:

\begin{proposition} \label{constr}
   Let $D$ be an effective purely exceptional divisor, such that for any
   $i,j \in I$, one has the inequality:
      $$(D+E_i+K_{\tilde{\mathcal{S}}})\cdot E_j + 2 \delta_i^j\leq 0$$
  where $\delta_i^j$ is Kronecker's symbol.
  Then there exists a function $f\in \frak{m}_{\mathcal{S},0}$ such
  that the exceptional part of $\mathrm{div}(f \circ \pi)$ is
  precisely $D$.
\end{proposition}

\textit{Proof:}
We use the following Grauert-Riemenschneider type vanishing theorem,
 proved by Laufer \cite{L 72} for analytic germs and by Ramanujam
 \cite{ R 72} for algebraic ones (see also B\u adescu \cite[4.1]{B 01}):
\medskip

 \textit{If $L$ is a divisor on $\tilde{\mathcal{S}}$ such that
 $L\cdot E_{j} \geq 
 K_{\tilde{\mathcal{S}}}\cdot E_{j},\:\forall\:j\in I$, then
 $H^{1}(\mathcal{O}_{\tilde{\mathcal{S}}}(L))=0$.  }
\medskip

 We apply the theorem to
$L=-D-E_i$, for any $i \in I$. Our hypothesis
implies that
$H^{1}(\mathcal{O}_{\tilde{\mathcal{S}}}(-D-E_i))=0$. Then, from the 
exact cohomology sequence associated to the short exact sequence of
sheaves
 $$ 0 \longrightarrow \mathcal{O}_{\tilde{\mathcal{S}}}(-D-E_i) \longrightarrow
 \mathcal{O}_{\tilde{\mathcal{S}}}(-D)
 \stackrel{\psi_{i}}{\longrightarrow}
 \mathcal{O}_{E_{i}}(-D)
 \longrightarrow 0$$
we deduce the \textit{surjectivity} of the restriction map
$$\psi_{i*}:H^{0}(\mathcal{O}_{\tilde{\mathcal{S}}}(-D)) \rightarrow
H^{0}(\mathcal{O}_{E_{i}}(-D))$$

By Serre duality on the irreducible
(possibly singular) curve $E_i$ (see Reid \cite[4.10]{R 97}), we get:
$$ h^{1}(\mathcal{O}_{E_i}(-D))
         = h^{0}(\mathcal{O}_{E_i}(K_{\tilde{\mathcal{S}}} +E_i +D))= 0$$
For the last equality we have used the hypothesis
   $(D+E_i+ K_{\tilde{\mathcal{S}}})\cdot E_i \leq -2 <0$, which shows
   that the line 
   bundle $\mathcal{O}_{E_i}(K_{\tilde{\mathcal{S}}} +E_i +D)$ cannot have a
   non-trivial section.

By applying the Riemann-Roch theorem and the adjunction formula for the
  irreducible curve $E_i$ of arithmetic genus $p_a(E_i)$ (see Reid
  \cite[4.11]{R 97}), we get:

$$
\begin{array}{ll}
   h^{0}(\mathcal{O}_{E_i}(-D)) &
      = h^{0}(\mathcal{O}_{E_i}(-D))- h^{1}(\mathcal{O}_{E_i}(-D))= \\
      & = \chi(\mathcal{O}_{E_i}(-D)) = \\
      & = 1 - p_a (E_i) - D\cdot E_i \geq \\
      & \geq 1 - p_a (E_i) + (K_{\tilde{\mathcal{S}}}+E_i)\cdot E_i +2= \\
      & = 1 -p_a (E_i) + 2 p_a (E_i) -2+2 = \\
      & =1 + p_a(E_i) >0
\end{array}
$$

This shows that there exists a non-identically zero section $s_{i}\in
H^{0}(\mathcal{O}_{E_{i}}(-D))$.
The surjectivity of
$\psi_{i*}$ implies that there exists $\sigma_{i} \in
H^{0}(\mathcal{O}_{\tilde{\mathcal{S}}}(-D))$ such that
$\psi_{i*}(\sigma_{i})=s_{i}$. As $(\mathcal{S},0)$ is normal, there
exists $f_{i} \in \frak{m}_{\mathcal{S},0}$ with
$\sigma_{i}=f_{i}\circ \pi$. If we write $D = \sum_{j \in I} a_j E_j$,
we see that $v_{E_i}(f_i)= a_i$ and $v_{E_j}(f_i)\geq a_j$, for all $j
\neq i$. We deduce that any generic linear combination
$f=\sum_{i\in I}\lambda_{i}f_{i}$ of the functions so constructed verifies:
  $$v_{E_i}(f) =a_i, \: \forall \: i \in I.$$
The proposition is proved.
  \hfill $\Box$

\medskip

\begin{remark} \textbf{a)} The proof follows the same line as the one
of Proposition 3.1 of \cite{CP 04} and 4.1 of \cite{CNP 04}, proved
by Caubel, N{\'e}methi and Popescu-Pampu. The
difference here is that we no longer deal with a \textit{good}
resolution of $(\mathcal{S},0)$ and we do not ask for a precise
knowledge of the topological type of the total transform of $f$.

\textbf{b)} As Morales informed us, Laufer \cite[3.1]{L 83} proved a
related statement when $\pi$ is the minimal resolution of
$(\mathcal{S},0)$: \emph{if $L$ is a line bundle on 
$\tilde{\mathcal{S}}_m$ with $L\cdot E_i \geq
2K_{\tilde{\mathcal{S}}_m}\cdot E_i, \: \forall 
\: i \in I$, then $L$ has no base points on $M$}. This implies
immediately, by taking $L=\mathcal{O}_{\tilde{\mathcal{S}}_m}(-D)$,
that the conclusion of 
Proposition \ref{constr} is true with the hypothesis $(D+2
K_{\tilde{\mathcal{S}}_m})\cdot E_i \leq 0, \: \forall \: i \in I$. As in
the sequel we work only with the minimal resolution, we could have
chosen to use Laufer's criterion. By looking at the way we use
Proposition \ref{constr} in the proof of Proposition \ref{separate},
one sees that this would have been enough in order to prove our main Theorem
\ref{extension}.  Nevertheless, we think that Proposition \ref{constr}
has independent interest.  Indeed, neither Laufer's hypothesis nor
ours implies the other one. In order to see it, consider first a
singularity whose minimal resolution has an irreducible exceptional
divisor $E$ which is smooth, of genus 2  and with $E^2=-1$: then
$D=4E$ verifies our hypothesis, but not Laufer's. Secondly, by
considering another singularity with $E= E_1 + E_2, 
\: p_a(E_1)=p_a(E_2)=0, \: E_1^2 =-4, \: E_2 ^2 =-2, \: E_1 \cdot E_2
=2$, one sees that $D =4E_1 + 4E_2$ verifies Laufer's hypothesis but
not ours. 
\end{remark}

\medskip

\section{The conditions (*) and (**)} \label{conditions}

From now on, we deal again with the minimal resolution $\pi_m$
of the germ $(\mathcal{S},0)$.

Inside the real vector space with basis  $(E_i)_{i\in I}$, we
consider the open half-spaces:
$$\{\sum_{i\in I} a_i E_i \: | \: a_i <a_j\},$$
for each $i\neq j$.  We call them
\textit{the fundamental half-spaces} of $\pi_m$.

We introduce  two conditions on the normal singularity
$(\mathcal{S},0)$:
\medskip

\textit{
(*) the intersection of $ \mathcal{L}^{\circ}(\pi_m)$ with each fundamental
half-space is non-empty;}

\textit{
(**) $E \in \mathcal{L}^{\circ}(\pi_m)$.}
\medskip

Notice that both conditions depend only on the intersection matrix of
$E$. In particular, they are purely topological.

\begin{lemma} \label{inclusion}
  The germs which verify condition (**) form a strict subset of those
  which verify condition (*).
\end{lemma}

\textit{Proof:} Let $(\mathcal{S},0)$ verify condition (**). Then, for
$n \in \mathbf{N}^*$ big enough (in fact for
$n >\displaystyle{\max_{i,j \in I} \{ \frac{E_i\cdot E_j}{|E\cdot
    E_i|}\}}$), one 
  sees from the definition of the strict Lipman semigroup that $nE+E_j
  \in \mathcal{L}^{\circ}(\pi_m), \: \forall
  \: j \in I$. But each fundamental half-space contains at least one
  of the divisors $nE+E_j$, which shows that $(\mathcal{S},0)$
  verifies condition (*).

Consider then any normal
   singularity $(\mathcal{S},0)$ whose minimal resolution has the same
   intersection matrix as a singularity of type $\mathbf{A}_n$ (that is, $E=
   \sum_{i=1}^n E_i$ with $E_i^2 =-2, \: \forall i\: \in \{1,...,n\},
   E_i\cdot E_{i+1}=1, \: \forall \: i\in \{1,...,n-1\}$ and $E_i\cdot
   E_j=0, \:
   \forall \: i,j \in \{1,...,n\}$ such that $|i-j| \notin \{0,1\}$), the
   components $E_i$ having otherwise arbitrary genera and
   singularities. Take $n \geq 3$. Then $E\cdot E_i =0, \: \forall \:
   i \in \{2,...,n-1\}$, which shows that $E \notin
   \mathcal{L}^{\circ}(\pi_m)$. That is, the germ $(\mathcal{S},0)$
   does not verify condition (**). 

Define
     $$ \alpha_k := nk -\frac{(k-1)k}{2}$$
   for any $k \in \{1,...,n\}$. Then it is immediate to see from the
   definitions
   that the divisors $D:=\sum_{k=1}^n \alpha_k  E_k$
   and $D':=  \sum_{k=1}^n \alpha_{n+1-k}E_k$ are in the strict Lipman
   semigroup and that each fundamental half-space contains exactly one of
   them. Therefore, the germ $(\mathcal{S},0)$ verifies condition (*). This
   shows that the inclusion stated in the Lemma is
   strict. \hfill $\Box$
\medskip

\begin{remark} \label{necessity}
    \textbf{a)} We choose to distinguish
   inside the class of singularities which
  verify condition (*) those which verify condition (**) for
   computational convenience, because this
  second condition is more readily verified on a given intersection
   matrix.

\textbf{b)} $\mathbf{A}_n$-type singularities are the only rational double
points which verify condition (*). This illustrates the difficulty of
dealing with $\mathbf{D}_n$-type singularities (see Pl{\'e}nat \cite{P 04}). 
\end{remark}

The motivation to introduce condition (*) comes from the following proposition:

\begin{proposition} \label{separate}
  Suppose that $(\mathcal{S},0)$ verifies condition (*). Then, for any
  pair of distinct
  indices $i,j \in I$, there exists a function $f \in
  \frak{m}_{\mathcal{S},0}$ such that $v_{E_i}(f) < v_{E_j}(f)$.
\end{proposition}

\textit{Proof:} Let $i,j \in I$ satisfy $i \neq j$. As
$(\mathcal{S},0)$ verifies condition (*), there exists $D =\sum_{l \in
  I}a_l E_l \in \mathcal{L}^{\circ}(\pi_m)$ such that $a_i <a_j$. Then, for
$n \in \mathbf{N}^*$ big enough, one has $(nD +E_k
+K_{\tilde{\mathcal{S}}}).E_l +2 \delta_k^l 
\leq 0, \: \forall \: k,l \in I$. By Proposition \ref{constr}, there
exists $f \in \frak{m}_{\mathcal{S},0}$ such that $\mathrm{div}(f
\circ \pi_m) = nD$, which shows that $v_{E_i}(f)=na_i < na_j
=v_{E_j}(f)$. \hfill $\Box$
\medskip

\begin{remark} \label{referee}
The referee suggests us the following alternative 
proof of the previous  proposition : let $D=\sum_{l \in  I}a_l
E_l\in\mathcal{L}^{\circ}(\pi_m)$ such  that $a_i <a_j$. Then
$\mathcal{O}_{\tilde{\mathcal{S}}_m}(-D)$ is ample  (see Lipman
\cite[10.4 and  proof 
of 12.1 (iii)]{L 69}). Thus there
exists $n \in \mathbf{N}^*$ such that
$\mathcal{O}_{\tilde{\mathcal{S}}_m}(-nD)$ 
is generated by its global sections.  Therefore there exists $f\in
\frak{m}_{\mathcal{S},0}$ such that the  exceptional part of
$\mathrm{div}(f \circ \pi)$ is precisely  $nD$. For such an $f$, we
have of course $v_{E_i}(f) < v_{E_j}(f)$. \end{remark}

\medskip

\section{The proof of the main theorem} \label{mainres}

Our main theorem is:

\begin{theorem} \label{extension}
  Suppose that $(\mathcal{S},0)$ verifies condition (*). Then the Nash map
  $\mathcal{N}$ is bijective.
\end{theorem}

\textit{Proof:} By combining Proposition \ref{distinguish} and
Proposition \ref{separate}, we deduce that $V(E_i) \nsubseteq
   V(E_j)$ for any $i\neq j$. As
  $(\mathcal{S},0)_{\infty}= \bigcup_{i \in I}V(E_i),$ 
we deduce that the schemes $(V(E_i))_{i \in I}$ are
precisely the irreducible components $(C_k)_{k\in K}$ of
$(\mathcal{S},0)_{\infty}$. As the Nash map $\mathcal{N}$ is
injective, this shows its surjectivity.

The Theorem is proved.
\hfill $\Box$
\medskip

Using Lemma \ref{inclusion}, we get as an immediate corollary:

\begin{corollary} \label{mthm2}
  Suppose that $(\mathcal{S},0)$ verifies condition (**). Then the Nash
  map $\mathcal{N}$  is bijective.
\end{corollary}

We denote by $\Gamma(E)$ the \textit{dual graph} of $E$, whose
vertices correspond bijectively to the components $(E_i)_{i \in
I}$, the vertex $E_i$ being weighted by $E_i^2$ and the vertices
$E_i$, $E_j$ being joined by $E_i.E_j$ vertices, for any $i \neq
j$. Let $\gamma(E_i)$ denote the number of edges which start from
the vertex $E_i$ (so, each loop based at the vertex $E_i$ counts for $2$).

The next proposition characterizes rational singularities among 
those which verify condition (**).

\begin{proposition} \label{rational}
 A singularity which verifies condition (**) is rational if and only
 if the following conditions are simultaneously verified:

  (i) \hspace{1.8mm} $\Gamma(E)$ is a tree;

  (ii) \hspace{0.5mm} $E_i \simeq \mathbf{P}^1, \: \forall \: i \in I$;

  (iii) $|E_i^2| > \gamma(E_i), \: \forall \: i \in I$.
\end{proposition}

\textit{Proof:} First of all, notice that conditions (i) and (ii)
imply that $\pi_m$ is a good resolution.

Suppose that $(\mathcal{S},0)$ is rational and
verifies condition (**). Then conditions (i) and (ii) are verified, as
general properties of rational singularities 
(see B\u{a}descu \cite[3.32.3]{B 01}). 
As $E.E_i= - |E_i^2| + \gamma(E_i)$, condition (iii) is also
verified.

Conversely, suppose that the conditions (i), (ii) and (iii) are
verified. By a result of Spivakovsky \cite[II]{S 90} (see also L{\^e} \cite[5.3]{L
97}), this shows that $(\mathcal{S},0)$ is a normal minimal
singularity, and in particular it is rational. Consult the references
above for the notion of minimal surface singularity, as well as
Koll{\'a}r \cite{K 85}, where this notion was introduced in arbitrary
dimensions. \hfill $\Box$

\medskip

\begin{remark}
  \textbf{a)} Normal minimal surface singularities are precisely those
  whose minimal good resolution verifies conditions (i), (ii)
  and (iii)': $|E_i^2| \geq \gamma(E_i)$,  $\forall \: i \in I$. The
  stronger condition (iii) is equivalent to the fact that the minimal
  resolution can be obtained by blowing-up once the origin (see
  Spivakovsky \cite[II]{S 90}, L{\^e} \cite[6.1]{L 97}): one says  that
  the singularity is 
  \textrm{superisolated}. Thus, Proposition \ref{rational} is equivalent
  to the fact that the rational surface singularities which verify
  condition (**) are precisely the superisolated minimal ones (a
  remark we owe to Lejeune-Jalabert).

  \textbf{b)} There are rational singularities which verify condition
  (*) but do
  not verify condition (**). Consider for example a germ of normal
  surface whose minimal resolution $\pi_m$ is good and has a reduced
  exceptional divisor $E$ with four components of genus $0$ such that
  $E_1^2 = E_2^2 =E_3^2 =-n \leq -5, \: E_4^2=-2, \: E_1 E_4 =E_2 E_4
  = E_3 E_4 =1, \: E_1 E_2 =E_2 E_3 =E_3 E_1 =0$. Then it is immediate
  to verify that the divisors:
$$ \begin{array}{c}(2n+1)(E_1 + E_2 +E_3)+ 4n E_4\\
      (2n^2 -2n +3)E_1 +3n(E_2 + E_3) +(n^2 +3n)E_4
   \end{array}$$
 as well as those obtained by permuting $E_1, E_2, E_3$ are in
      $\mathcal{L}^{\circ}(\pi_m)$ and that each fundamental
      half-space contains at least one of them. So, the germ verifies
      condition (*). But, as $E.E_4 >0$,  the germ does not verify condition
      (**). Moreover, the dual graph $\Gamma(E)$ is a subgraph of the
      dual graph associated to the resolution of a plane curve
      singularity (attach
      $n-1$ vertices to $E_1$, $n-2$ vertices to $E_2$ and $n-3$
      vertices to $E_3$, all of them weighted by $-1$). Using L{\^e}
      \cite[4.8]{L 97} (see also Spivakovsky \cite[II]{S 90}), we see that
      the singularity is 
      sandwiched, and in particular it is rational.
\end{remark}

As condition (iii) in the Proposition \ref{rational} is equivalent with
condition (**), we see
immediately that any (abstract) graph is the dual graph of the minimal
resolution of a singularity which verifies condition (**), once the
weights of the vertices are negative enough. Then, if one of the
conditions (i) or (ii) is not satisfied, we are in presence of a
non-rational singularity. This shows:

\begin{corollary} \label{infinity}
  There exists an infinity of pairwise topologically distinct
  normal non-rational surface singularities which verify condition
  (**), and consequently for which the Nash map is bijective.
\end{corollary}

\medskip

\end{document}